\documentclass[11pt, 14paper,reqno]{amsart}
\usepackage{amsmath,amsfonts,amssymb}
\usepackage{algorithm}
\usepackage{algpseudocode}
\usepackage{bigints}
\usepackage{longtable}
\usepackage{color}
\makeatletter
\newcommand{\least}{\let\CS=\@currsize\renewcommand{\baselinestretch}{.9}\tiny\CS}
\renewcommand\baselinestretch{1.1} 
\theoremstyle{plain}
\newtheorem{theorem}{Theorem}

\theoremstyle{proof}
\theoremstyle{definition}

\theoremstyle{remark}
\newtheorem{remark}{Remark}
\theoremstyle{lamma}

\numberwithin{equation}{section}
\numberwithin{lemma}{section}
\numberwithin{theorem}{section}
\numberwithin{remark}{section}
\numberwithin{prop}{section}
\numberwithin{corollary}{section}

\usepackage{amsmath}
\usepackage{amsfonts}   
\usepackage{amssymb}
\usepackage{amssymb, amsmath, amsthm}
\usepackage[breaklinks]{hyperref}

\theoremstyle{thmrm} 
\newtheorem{exa}{Example}

\numberwithin{conjecture}{section}

\begin{document}
\title[Computation of Jacobi sums and cyclotomic numbers with reduced complexity]{\large{Computation of Jacobi sums and cyclotomic numbers with reduced complexity}}
\author{Md Helal Ahmed and Jagmohan Tanti}
\address{Md Helal Ahmed @ Department of Mathematics, Central University of Jharkhand, Ranchi-835205, India}
\email{ahmed.helal@cuj.ac.in}

\address{Jagmohan Tanti @ Department of Mathematics, Central University of Jharkhand, Ranchi-835205, India}
\email{jagmohan.t@gmail.com}

\maketitle
\begin{abstract}
Jacobi sums and cyclotomic numbers are the important objects in number theory. The determination of all the Jacobi sums and cyclotomic numbers of order $e$ are merely intricate to compute. This paper presents the lesser numbers of Jacobi sums and cyclotomic numbers which are enough for the determination all Jacobi sums and the cyclotomic numbers of a particular order.
\end{abstract}
\section{Introduction}\label{sec1}
 For $e\in \mathbb{Z}^{+}\geq 2$, in a finite field a Jacobi sum of order $e$ mainly depends on two parameters. Therefore, these values could be naturally assembled into a matrix of order $e$. As illustrated in \cite{Ireland1}, Jacobi sums could be used for estimating the number of integral solutions to congruences such as $x^{3}+y^{3}\equiv 1 \pmod p$. These estimates play a key role in the development of the Weils conjecture \cite{Weil1}. Jacobi sums were also utilized by Adleman, Pomerance, Rumely \cite{Adleman1} for primality testing.

Let $e\geq 2$ be an integer, $p$ a prime, $q=p^{r}, r \in\mathbb{Z}^{+}$ and $q\equiv 1 \pmod{e}$. One writes $q=ek+1$ for some positive integer $k$. Let $\mathbb{F}_{q}$ be a finite field of $q$ elements and $\gamma$ a generator of the cyclic group $\mathbb{F}_{q}^{*}$. For a primitive $e$-th root $\zeta_e$ of unity, define a multiplicative character $\chi_{e}$ of order $e$ on $\mathbb{F}^{*}_{q}$ by $\chi_{e}(\gamma)=\zeta_e$. We now extent $\chi_{e}$ to a map from $\mathbb{F}_q$ to $\mathbb{Q}({\zeta_e})$ by taking $\chi_{e}(0)=0$. For $0\leq i, j\leq e-1$, the Jacobi sums of order $e$ is defined by 
$$
J_e(i,j)= \sum_{v\in \mathbb{F}_q} \chi_{e}^i(v) \chi_{e}^j(v+1).
$$
\par 
For $0\leq a, b\leq e-1$, the cyclotomic numbers $(a, b)_e$ of order $e$ is defined as follows:
\begin{align*}
(a,b)_e:& =\#\{v\in\mathbb{F}_q|\chi_e(v)=\zeta_e^a,\,\,\chi_e(v+1)=\zeta_e^b\} \\ & =
\#\{v\in\mathbb{F}_q\setminus \{0,-1\}\mid {\rm ind}_{\gamma}v\equiv a \pmod e ,  {\ \rm ind}_{\gamma}(v+1)\equiv b \pmod e\}.
\end{align*}
The Jacobi sums $J_e(i,j)$ and the cyclotomic numbers $(a,b)_e$ are well connected  by the following relations \cite{Berndt1,Shirolkar1}:
\begin{equation} \label{01}
 \sum_a\sum_b(a,b)_e\zeta_e^{ai+bj}=J_e(i,j),
\end{equation} 
and
\begin{equation} \label{00}
 \sum_i\sum_j\zeta_e^{-(ai+bj)}J_e(i,j)=e^{2}(a,b)_e.
\end{equation}
For more about the properties of Jacobi sums and cyclotomic numbers, one may refer \cite{Berndt1,Shirolkar1,Acharya1}. 
\par
Cyclotomic numbers are one of the most important objects in number theory. These numbers have been extensively used in cryptography, coding theory and other branches of information theory. Thus determination of cyclotomic numbers, so called cyclotomic number problem, of different order is one of basic problems in number theory. Complete solutions to cyclotomic number problem for $e$ =  $2-6$, $7$, $8$, $9$, $10$, $11$, $12$, $14$, $15$, $16$, $18$, $20$, $22$, $l$, $2l$, $l^{2}$, $2l^{2}$ with $l$ an odd prime have been investigated by many authors \cite{Shirolkar1,Acharya1,Katre3}.
\par
The problem of determining Jacobi sums and cyclotomic numbers has been treated by many authors since the age of Gauss. Jacobi sums and cyclotomic numbers are interrelated which is shown in equations (\ref{01}) and (\ref{00}). The problem of cyclotomic numbers of prime order $l$ in the finite field $\mathbb{F}_{q}$, $q=p^{r}$, $p \equiv 1 \pmod l$ has been investigated by Katre ans Rajwade \cite{Katre3} in terms of the arithmetic characterization of Jacobi sums of prime order $l$. Later, Acharya and Katre \cite{Acharya1} induced some additional condition and solve the cyclotomic problem twice of a prime in terms of the coefficients of Jacobi sums of orders $l$ and $2l$. Recently, Shirolkar and Katre \cite{Shirolkar1} solved the cyclotomic problem of order $l^{2}$ in terms of the coefficients of Jacobi sums of orders $l$ and $l^{2}$. 
\par
The question of determining all Jacobi sums of orders $e$ in terms of lesser number of Jacobi sums has also been taken by many authors such as Parnami, Agrawal and Rajwade \cite{Parnami2}. Similar kind of arguments for cyclotomic numbers can also be given. In this paper, we give a list of fewer Jacobi sums (respectively cyclotomic numbers) of order $e$ which qualifies to determine all Jacobi sums (respectively cyclotomic numbers) of order $e$. We exhibit an explicit case with $e=3$.  

\section{Main Results}
\subsection{Computation of Cyclotomic Numbers of order $e$} \label{sec6}
The following result is a sort of optimization of the cyclotomic number problem. 
\begin{theorem}
The number of cyclotomic numbers $(a,b)_{e}$ required for the determination of all the cyclotomic numbers of order $e$ is equal to $e+(e-1)(e-2)/6$, if $6|(e-1)(e-2)$, otherwise $e+\Big \lceil (e-1)(e-2)/6 \Big\rceil +1$.
\end{theorem}
\begin{proof}
Recall from the properties of cyclotomic numbers \cite{Berndt1,Shirolkar1,Acharya1}, if $k$ is even or $q=2^r$, then 
\begin{equation} \label{equ8}
(a,b)_{e}=(b,a)_{e}=(a-b,-b)_{e}=(b-a,-a)_{e}=(-a,b-a)_{e}=(-b,a-b)_{e}
\end{equation}
otherwise
\begin{align} \label{equ9}
\nonumber &(a,b)_{e}=(b+\frac{e}{2},a+\frac{e}{2})_{e}=(\frac{e}{2}+a-b,-b)_{e}=(\frac{e}{2}+b-a,\frac{e}{2}-a)_{e} \\ &=(-a,b-a)_{e}  =(\frac{e}{2}-b,a-b)_{e}.
\end{align}
Thus by (\ref{equ8}) and (\ref{equ9}) partition cyclotomic numbers $(a,b)_{e}$ of order $e$ into group of classes. The expression (\ref{equ8}) splits the problem into two cases: 

\textbf{Case 1:} $2|k$ $/$ $q=2^r$ and $3 \nmid e$: In this case, $(0,0)_{e}$ form singleton class, $(-a,0)_{e}$, $(a,a)_{e}$, $(0,-a)_{e}$ form classes of three elements where $1\leq a\leq$ ${e-1 \pmod {e}}$ and rest $(e^{2}-3e+2)$ of the cyclotomic numbers form classes of six elements. Therefore the total number of distinct cyclotomic number classes is equal to $1+(e-1)+(e^{2}-3e+2)/6$ $=$ $e+(e-1)(e-2)/6$.   

\textbf{Case 2:} $2|k$ $/$ $q=2^r$ and $3|e$ with $e=3x$ for some $x \in \mathbb{Z}^{+}$: In this situation, $(0,0)_{e}$ form singleton class, $(-a,0)_{e}$, $(a,a)_{e}$, $(0,-a)_{e}$ form classes of three elements where $1\leq a\leq e-1 \pmod {e}$, $(x,2x)_{e}, (2x,x)_{e}$ which are grouped into classes of two elements and rest $(e^{2}-3e)$ of the cyclotomic numbers form classes of six elements. Therefore the total number of distinct cyclotomic number classes is equal to $1+1+(e-1)+(e^{2}-3e)/6$ $=$ $e+\Big \lceil (e-1)(e-2)/6 \Big\rceil +1$.

Again the expression (\ref{equ9}) splits the problem into two cases: 

\textbf{Case 3:} $2\nmid k$ $\&$ $q\neq 2^r$ and $3\nmid e$: In this case, $(0,\frac{e}{2})_{e}$ form singleton class, $(0,a)_{e}$, $(a+\frac{e}{2},\frac{e}{2})_{e}$, $(\frac{e}{2}-a,-a)_{e}$ form classes of three elements where $1\leq a\leq {e}-1 \pmod {e}$ and rest $(e^{2}-3e+2)$ of the cyclotomic numbers form classes of six elements. Therefore the total number of distinct cyclotomic number classes is equal to $1+(e-1)+(e^{2}-3e+2)/6$ $=$ $e+(e-1)(e-2)/6$. 

\textbf{Case 4:} $2\nmid k$ $\&$ $q\neq 2^r$ and $3|e$: (it implies that $e=6x$ for some $x \in \mathbb{Z}^{+}$). Here cyclotomic numbers $(0,\frac{e}{2})_{e}$ form singleton class, $(0,a)_{e}$, $(a+\frac{e}{2},\frac{e}{2})_{e}$, $(\frac{e}{2}-a,-a)_{e}$ form classes of three elements where $1\leq a\leq e-1 \pmod {e}$, $(2x,x)_{e}, (x+\frac{e}{2},2x+\frac{e}{2})_{e}$ which are grouped into classes of two elements and rest $(e^{2}-3e)$ of the cyclotomic numbers form classes of six elements. Thus the total number of distinct cyclotomic numbers is equal to $1+1+(e-1)+(e^{2}-3e)/6$ $=$ $e+\Big \lceil (e-1)(e-2)/6 \Big\rceil +1$.
\par
Therefore, the number of cyclotomic numbers required for the determination of all the cyclotomic numbers for $e\geq 2$, is equal to $e+(e-1)(e-2)/6$, if $6|(e-1)(e-2)$, otherwise $e+\Big \lceil (e-1)(e-2)/6 \Big\rceil +1$. 
\end{proof}
\begin{remark}
For the determination all the Jacobi sums $J_{e}(i,j)$ of order $e$, it is sufficient to determine exactly $e+(e-1)(e-2)/6$ numbers of cyclotomic numbers of order $e$, if $(e-1)(e-2)$ completely divisible by $6$, otherwise $e+\Big \lceil (e-1)(e-2)/6 \Big\rceil +1$.
\end{remark}
\subsection{Computation of Jacobi sums of order $e$}\label{sec4}
\begin{theorem}
The number of Jacobi sums $J_{e}(i,j)$ required for the determination of all the Jacobi sums of order $e$ is equal to $e+(e-1)(e-2)/6$, if $6|(e-1)(e-2)$, otherwise $e+\Big \lceil (e-1)(e-2)/6 \Big\rceil +1$.
\end{theorem}
\begin{proof}
From the properties of Jacobi sums \cite{Berndt1,Shirolkar1,Acharya1}, if $k$ is even or $q=2^r$, then
\begin{equation} \label{equation3}
J_{e}(i,j)=J_{e}(j,i)=J_{e}(-i-j,i)=J_{e}(i,-i-j)=J_{e}(j,-i-j)=J_{e}(-i-j,j),
\end{equation}
otherwise
\begin{align*}
&(-1)^{ik}J_{e}(i,j)=(-1)^{jk}J_{e}(j,i)=(-1)^{jk}J_{e}(-i-j,i)=(-1)^{jk}J_{e}(i,-i-j)\\ & =(-1)^{ik}J_{e}(j,-i-j)=(-1)^{ik}J_{e}(-i-j,j).
\end{align*}
These above expressions partition the Jacobi sums $J_{e}(i,j)$ of order $e$ into group of classes. It is clear that Jacobi sums matrix is always symmetric and if $k$ is even or $q=2^r$ or $k$ odd, the entries of the Jacobi sums matrix differ atmost by sign. Let us consider the expression (\ref{equation3}). It split the problem into three cases:

\textbf{Case 1:} $i=j:$ In this case, Jacobi sums of order $e$ partition into classes of singleton and three elements. If $3i\equiv 0 \pmod e$ form singleton class i.e. $J_{e}(i,i)$, otherwise form classes of three elements i.e.$J_{e}(i,i)$ $=$ $J_{e}(i,-2i)$ $=$ $J_{e}(-2i,i)$. 

\textbf{Case 2:} $i=0$ and $i\neq j$: In this situation, Jacobi sums of order $e$ partition into classes of three and six elements. If $2j\equiv 0 \pmod e$ form classes of three elements i.e. $J_{e}(0,j)$ $=$ $J_{e}(j,0)$ $=$ $J_{e}(-j,j)$, otherwise form classes of six elements i.e.$J_{e}(0,j)$ $=$ $J_{e}(j,0)$ $=$ $J_{e}(-j,j)$ $=$ $J_{e}(j,-j)$ $=$ $J_{e}(-j,0)$ $=$ $J_{e}(0,-j)$.

\textbf{Case 3:} $i\neq 0$ and $i\neq j$: Again in this case, Jacobi sums of order $e$ partition into classes of three and six elements. If $i+2j\equiv 0 \pmod e$ form classes of three elements i.e. $J_{e}(i,j)$ $=$ $J_{e}(j,i)$ $=$ $J_{e}(-i-j,j)$, otherwise form classes of six elements i.e.$J_{e}(i,j)$ $=$ $J_{e}(j,i)$ $=$ $J_{e}(-i-j,j)$ $=$ $J_{e}(j,-i-j)$ $=$ $J_{e}(-i-j,i)$ $=$ $J_{e}(i,-i-j)$.
\par
Thus, it is clear by above cases, for the determination of all the Jacobi sums for $e\geq 2$, it is enough to determine $e+(e-1)(e-2)/6$ number of Jacobi sums, if $6|(e-1)(e-2)$, otherwise $e+\Big \lceil (e-1)(e-2)/6 \Big\rceil +1$.   
\end{proof}
\begin{remark}
For the determination all the cyclotomic numbers $(a,b)_{e}$ of order $e$, it is sufficient to determine exactly $e+(e-1)(e-2)/6$ numbers of Jacobi sums $J_{e}(i,j)$ of order $e$, if $6|(e-1)(e-2)$, otherwise $e+\Big \lceil (e-1)(e-2)/6 \Big\rceil +1$.
\end{remark}
Here we give an example for optimize of computation for determination of Jacobi sums and cyclotomic numbers respectively of order $3$.
\begin{exa}
Let us choose a number $q$ that satisfy $q=p^{r}\equiv 1 \pmod 3$, where $p$ prime and $r\in \mathbb{Z}^{+}$. Assume that $q=7$. The number of cyclotomic numbers required for the determination all the cyclotomic numbers of order $3$ is shows in table \ref{table1}.
{\tiny{
\begin{table}[h!] \label{table1}
\noindent\makebox[-5mm]{
\begin{tabular}{|p{0.5cm}| p{0.6cm} p{0.6cm} p{0.6cm} |}\hline(a,b) &\multicolumn{3}{c|}{b}\\ \hline
a & 0 & 1 & 2
\\ \hline
0 &(0,0)&(0,1)&(0,2)
\\ \hline
1 &(0,1)&(0,2)&(1,2)
\\ \hline
2 &(0,2)&(1,2)&(0,1)
\\ \hline
\end{tabular}}
\vspace*{2mm}
\caption{}
\end{table}}}

For the determination of all the Jacobi sums of order $3$ over $\mathbb{F}_{7}$, we need to determine only four cyclotomic numbers out of nine cyclotomic numbers, which are $(0,0)_{3}$, $(0,1)_{3}$, $(0,2)_{3}$ and $(1,2)_{3}$. For the determination of above cyclotomic numbers, apply the formula \cite{Dickson1}, we get $(0,0)_{3}=0$, $(0,1)_{3}=0$, $(0,2)_{3}=1$ and $(1,2)_{3}=1$ w.r.t. generator $3$ of $\mathbb{F}_{7}^{*}$, which is shown in matrix $A$. 
\begin{center}
	
	{\small$A
	=\left[
	\begin{array}{rrr}
    0 & 0 & 1 
	\\
	0 & 1 & 1  
	\\
	1 & 1 & 0 
	\\
	\end{array}\right] 
	$ }
\end{center}   
Apply the relation (\ref{01}) for the determination of $J_{3}(i,j)$, where $0\leq i,j\leq 2$,
\begin{equation*}
J_{3}(i,j)= \zeta_{3}^{2j}+\zeta_{3}^{i+j}+\zeta_{3}^{i+2j}+\zeta_{3}^{2i}+\zeta_{3}^{2i+j}.
\end{equation*}
From \cite{Parnami2}, we have $J_{3}(0,0)=q-2$, $J_{3}(0,1)=-1$, $J_{3}(0,2)=-1$, $J_{3}(1,0)=-1$, $J_{3}(1,2)=-1$, $J_{3}(2,0)=-1$, $J_{3}(2,1)=-1$. Now we need to determine only $J_{3}(1,1)$ and $J_{3}(2,2)$.
\begin{equation*}
J_{3}(1,1)= 3\zeta_{3}^{2}+2\zeta_{3}^{3}.
\end{equation*}
For further simplification, employ the cyclotomic polynomial of order $3$, which is $\zeta_{3}^{2}+\zeta_{3}+1$.  Finally, we get \begin{equation*}
J_{3}(1,1)=-3\zeta_{3}-1.
\end{equation*}
Similarly, 
\begin{equation*}
J_{3}(2,2)=3\zeta_{3}+1.
\end{equation*} 
Conversely, the number of Jacobi sums required for the determination all the Jacobi sums of order $3$ is shows in table \ref{table2}.
{\tiny{
\begin{table}[h!] 
\noindent\makebox[-5mm]{
\begin{tabular}{|p{0.5cm}| p{0.6cm} p{0.6cm} p{0.6cm} |}\hline(a,b) &\multicolumn{3}{c|}{b}\\ \hline
a & 0 & 1 & 2
\\ \hline
0 &(0,0)&(0,1)&(0,1)
\\ \hline
1 &(0,1)&(1,1)&(0,1)
\\ \hline
2 &(0,1)&(0,1)&(2,2)
\\ \hline
\end{tabular}}
\vspace*{2mm}
\caption{} \label{table2}
\end{table}}}

\indent Now, for the determination of all the cyclotomic numbers of order $3$ over $\mathbb{F}_{7}$ in terms of Jacobi sums of order $3$, we need to determine only four Jacobi sums of order $3$ out of nine Jacobi sums of order $3$, these are $J_{3}(0,0)$, $J_{3}(0,1)$, $J_{3}(1,1)$ and $J_{3}(2,2)$. Now, for the determination of $(a,b)_{3}$, where $0\leq a,b\leq 2$, apply the relation (\ref{00}). On substituting the values of $J_{3}(0,0)$, $J_{3}(0,1)$, $J_{3}(1,1)$ and $J_{3}(2,2)$, we get $(0,0)_{3}=0$, $(0,1)_{3}=0$, $(0,2)_{3}=1$, $(1,0)_{3}=0$, $(1,1)_{3}=1$, $(1,2)_{3}=1$, $(2,0)_{3}=1$, $(2,1)_{3}=1$, $(2,2)_{3}=0$.
\end{exa}
\section{Conclusion}\label{sec5}
We expect that the calculated number of cyclotomic numbers and Jacobi sums are the minimal for the determination cyclotomic numbers and Jacobi sums, as lesser than this has not been seen in the literature.
\vskip2mm
\noindent{\bf Acknowledgement}
The authors are thankful to the Central University of Jharkhand, Ranchi, India for providing the research support.

\end{document}